\documentclass{birkjour}
\usepackage{amsmath,amssymb,amsthm}

\theoremstyle{definition}
\newtheorem{theorem}{Theorem} 
\newtheorem{remark}{Remark}[section] 
\newtheorem{definition}{Definition}[section]

\DeclareMathOperator{\Li}{Li}

\begin{document}

\title[An identity for the Bernoulli numbers]
 {An identity involving Bernoulli numbers and the Stirling numbers of the second kind}

\author[S. K. Jha]{Sumit Kumar Jha}

\address{%
International Institute of Information Technology\\
Hyderabad-500 032, India}

\email{kumarjha.sumit@research.iiit.ac.in}





\begin{abstract}
Let $B_{n}$ denote the Bernoulli numbers, and $S(n,k)$ denote the Stirling numbers of the second kind. We prove the following identity
 $$ B_{m+n}=\sum_{\substack{0\leq k \leq n \\ 0\leq l \leq m}}\frac{(-1)^{k+l}\,k!\, l!\, S(n,k)\,S(m,l)}{(k+l+1)\,\binom{k+l}{l}}. $$
To the best of our knowledge, the identity is new.
\end{abstract}

\maketitle

\section{Introduction} 

\begin{definition}
The \emph{Bernoulli numbers} $B_{n}$ can be defined by the 
following generating function:
\begin{equation*}
    \frac{t}{e^{t}-1}=\sum_{n\geq 0}\frac{B_{n}t^{n}}{n!},
\end{equation*}
where $|t|<2\pi$.
\end{definition}
\begin{definition}
The \emph{Stirling number of the second kind}, denoted by $S(n,m)$, is the number of ways of partitioning a set of $n$ elements into $m$ nonempty sets.
\end{definition}
The following formula expresses the Bernoulli numbers explicitly in terms of the Stirling numbers of the second 
kind \cite{Gould,Qi1}:
\begin{align}
   B_{n}&=\sum_{k=0}^{n}\frac{(-1)^{k}\, k! \, S(n,k)}{k+1}. \label{firstBern} 
\end{align}
In the following section, we prove a new identity for the Bernoulli numbers in terms of Stirling numbers of the second kind, of which the above formula is a special case.
\section{Main result}
Our main result is the following.
\begin{theorem}
For all non-negative integers $m,n$ we have
\begin{equation*}
B_{m+n}=\sum_{\substack{0\leq k \leq n \\ 0\leq l \leq m}}\frac{(-1)^{k+l}\, k!\, l!\, S(n,k)\,S(m,l)}{(k+l+1)\,\binom{k+l}{l}}.
\end{equation*}
\end{theorem}
\begin{remark}
Letting $m=0$ in above equation gives us equation \eqref{firstBern}.
\end{remark}
\begin{proof}
We start with the following integral from \cite{Wolfram}
\begin{equation}
(\alpha+\beta)\zeta(\alpha+\beta+1)=\int_{0}^{\infty} \frac{\Li_{\alpha}(-1/t)\, \Li_{\beta}(-t)}{t}\, dt,
\end{equation}
where $\zeta(\cdot)$ is the Riemann zeta function, and $\Li_{\alpha}(t)$ is the polylogarithm function. \par
Letting $\alpha=-m$, and $\beta=-n$ (non-negative integers) in the preceding equation, we get
\begin{equation*}
-(m+n)\zeta(1-m-n)=\int_{0}^{\infty} \frac{\Li_{-m}(-1/t)\, \Li_{-n}(-t)}{t}\, dt.
\end{equation*}
The following representation from the note \cite{Landsburg}
\begin{equation}
\label{stirling}
    \Li_{-n}(-t)=\sum_{k=0}^nk!\, S(n,k)\left(\frac{1}{1+t}\right)^{k+1}(-t)^{k}
\end{equation}
allows us to evaluate the integral as
\begin{align*}
\int_{0}^{\infty} \frac{\Li_{-m}(-1/t)\, \Li_{-n}(-t)}{t}\, dt &= 
\int_{0}^{\infty}\sum_{\substack{0\leq k \leq n \\ 0\leq l \leq m}}\frac{(-1)^{k+l}\,k!\,l! \, S(n,k)\, S(m,l)\, t^{k}}{(1+t)^{k+l+2}}\, dt\\
&=\sum_{\substack{0\leq k \leq n \\ 0\leq l \leq m}}(-1)^{k+l}\,k!\,l! \, S(n,k)\, S(m,l)\int_{0}^{\infty}\frac{t^{k}}{(1+t)^{k+l+2}}\, dt\\
&=\sum_{\substack{0\leq k \leq n \\ 0\leq l \leq m}}(-1)^{k+l}\,k!\,l! \, S(n,k)\, S(m,l)\frac{\Gamma(k+1)\Gamma(l+1)}{\Gamma(k+l+2)}.
\end{align*}
Here $\Gamma(\cdot)$ is the gamma function. This completes the proof after noting the fact \cite{Ahlfors} that $$-(m+n)\cdot \zeta(1-m-n)=B_{m+n}.$$ 
\end{proof}

\makeatletter
\renewcommand{\@biblabel}[1]{[#1]\hfill}
\makeatother

\end{document}